\documentclass[12pt,a4paper]{article}

\usepackage{amssymb}

\newtheorem{Proposition}{Proposition}
\newtheorem{Theorem}[Proposition]{Theorem}
\newtheorem{Corollary}[Proposition]{Corollary}
\newtheorem{Lemma}[Proposition]{Lemma}

\newcommand{\C}{{\mathbb C}}
\newcommand{\N}{\mathbb N}
\newcommand{\Z}{\mathbb Z}
\newcommand{\R}{\mathbb R}
\newcommand{\Qp}{\mathbb Q_p} 

\newcommand {\F}{\mbox{${\cal F}$}}

\newcommand{\proof}{\noindent {\bf Proof: }}
\newcommand {\qed}{\hfill $\square$}

\begin{document}

\title{On Extensions of generalized Steinberg Representations}
\author{Sascha Orlik}
\date{}
\maketitle

\begin{abstract}
Let $F$ be a local non-archimedean field and let $G$ be the group of
$F$-valued points of a reductive algebraic group over $F.$
In this paper we compute the Ext-groups of generalized Steinberg representations
in the category of smooth $G$-representations with coefficients in a certain self-injective ring. 

\end{abstract}

\section{Introduction}
The origin of the problem we treat here is the computation of the
\'etale cohomology of $p$-adic period domains with finite coefficients. 
In \cite{O} the computation yields a filtration of smooth
representations of a $p$-adic Lie group on the cohomology groups,
which is induced by a certain spectral sequence. 
A natural problem which arises in this context is to show that this filtration splits canonically. The graded pieces 
of the filtration are essentially generalized Steinberg representations. A natural task is 
therefore to study the extensions of these  representations.

Let $F$ be a local non-archimedean field and let $G$ be the group of
$F$-valued points of a fixed reductive algebraic group over $F.$ The
field $F$ induces a natural topology on $G$ providing it with the structure
of a locally profinite group.
The aim of this paper is to determine the Ext-groups of 
generalized Steinberg representations 
in the category of smooth $G$-representations with coefficients in 
a self-injective ring $R$. We refer to the next chapter
for the precise conditions we impose on $R.$ An important example of
such a ring is given by a field of characteristic zero. 
One crucial assumption is that the pro-order of $G$ is invertible
in $R.$ In \cite{V1} it is shown that this condition is sufficient  
for the existence of a (left-invariant) normalized Haar measure on $G.$
Using this Haar measure and the self-injectivity of $R$ ensures all the well-known properties and techniques
in representation- and cohomology theory of a $p$-adic reductive group, e.g. Frobenius reciprocity, exactness of the fixed point functor for a compact
open subgroup of $G$  etc.,  as in the classical case where $R =\C.$ 
In particular we have enough injective and projective objects in the
category of smooth $G$-representations.

The generalized Steinberg representations
are parametrized by the subsets of a relative $\Qp$-root basis $\Delta$ of $G.$ 
For any subset $I\subset \Delta,$ let $P_I\subset  G$ the corresponding standard-parabolic subgroup of $G.$ 
Let $i^G_{P_I}=C^\infty(P_I\backslash G)$ be the 
$G$-representation consisting of locally constant functions on
$P_I\backslash G$ with values in $R.$
If $J\supset I$ is another subset, then there is  a natural injection $i^G_{P_J}
\hookrightarrow i^G_{P_I}.$
The generalized Steinberg representation with respect to $I\subset
\Delta$ is the quotient
$$v^G_{P_I}= i^G_{P_I}/\sum_{I\subset J\subset \Delta \atop I\neq J} i^G_{P_J}.$$ In the
case $I=\emptyset$ we just get the ordinary Steinberg representation.
In the case $R=\C$ it is known that the representations $v^G_{P_J},$
for $J\supset I,$ are precisely the
irreducible subquotients of $i^G_{P_I}.$
Our main result is formulated in the following theorem.

\begin{Theorem}  Let $G$ be semi-simple. Let $I,J \subset \Delta.$ Then
$$Ext_G^i(v^G_{P_I},v^G_{P_J})=\left\{ \begin{array}{r@{\quad:\quad}l} R & i = |I\cup J| -  |I\cap J| \\ 0 & otherwise \end{array} \right. .$$
\end{Theorem}
Note that in the case where $I$ or $J$ is the empty set, i.e., $v^G_{P_I}$ or $v^G_{P_J}$ is the trivial representation
and $R$ is the field of complex numbers,  this computation has been carried out by
Casselman \cite{Ca1}, \cite{Ca2} resp. Borel and Wallach \cite{BW}. If on the other extreme $I=\Delta$ or $J=\Delta$, the Ext-groups have been computed by Schneider and Stuhler \cite{SS}.

If $G$ is not necessarily semi-simple then we also have a
contribution of the center $Z(G)$ of $G.$  
By using a Hochschild-Serre argument we conclude from Theorem 1:
\begin{Corollary}
Let $G$ be reductive with center $Z(G)$ of $\Qp$-rank $d.$  Let $I,J \subset \Delta.$ Then we have
$$Ext_G^i(v^G_{P_I},v^G_{P_J})=\left\{ \begin{array}{r@{\quad:\quad}l} R^{d \choose j }  & i = |I\cup J| -  |I\cap J| +j \\ 0 & otherwise \end{array} \right. $$
\end{Corollary}

During my computations I was informed by J.-F. Dat that
he was also able to prove Theorem 1. His
proof \cite{D} is totally different from ours. It
is based on intertwining operators and Bernstein's second adjunction
formula. In addition to the fact that $R$ need not to be
self-injective, his proof has the advantage of producing
the extensions of generalized Steinberg representations explicitly.  

Our proof of Theorem 1 is quite natural. One uses certain resolutions of the representations $v^G_{P_I}$ in terms of the induced representations
$i^G_{P_K},$ where $K\supset I.$  By a spectral sequence argument, the proof reduces to  the computation of the groups 
$Ext_G^\ast(i^G_{P_I},i^G_{P_J}),$ for  $I,J\subset \Delta.$
This is done by Frobenius reciprocity and a description of the Jacquet
modules for these kind of representations.
The latter  has been considered in \cite{Ca3} in the case $R =\C.$ 
It holds more generally in our situation.

I am grateful to J.-F. Dat  for his numerous remarks on this paper. He
explained to me how to genaralize my proof from the case
$R=\C$ to the case of a certain self-injective ring.
I would like to thank the IHES and J.-F. Dat  for the invitation in
June 2003. I wish to thank A. Huber and M. Rapoport for helpful remarks.
I also thank T. Wedhorn and P. Schneider for their comments on a first version of this paper.
Finally, I would like to thank C. Kaiser for pointing out to me
Corollary 18 as a consequence of the results above.

\section{Notations}
Let $p$ be a prime number and let $F$ be a local non-archimedean
field. We suppose that the residue field of $F$ has order $q=p^r,
r\geq 0.$ Let $val:F\rightarrow \Z$
be the discrete valuation taking a fixed uniformizer $\varpi_F\in F$ to $1\in \Z.$
Denote by 
$| \;\; |_{\R} : F \rightarrow \R$ the corresponding
normalized $p$-adic norm with values in $\R$.

Let ${\bf G}$ be a reductive algebraic group over $F.$
Fix a maximal $F$-split torus ${\bf S}$  and a minimal $F$-parabolic subgroup ${\bf P}$ in ${\bf G}$ containing ${\bf S}.$ 
Let ${\bf M} = Z({\bf S})$ be the centralizer of ${\bf S}$
in ${\bf G}$, which is a Levi subgroup of ${\bf P}.$
Denote by ${\bf U}$ the unipotent radical of ${\bf P}.$
Let 
$$\Phi\supset\Phi^+\supset\Delta=\{\alpha_1,\ldots,\alpha_n\}$$
 be the corresponding subsets of relative $F$-roots, $F$-positive roots, $F$-simple roots. 
In the following, we call them for simplicity just roots instead of  relative $F$-roots.
For a subset $I\subset \Delta,$ we let ${\bf P}_I \subset {\bf G}$
 be the standard parabolic subgroup defined over $F$ such
 that $\Delta \setminus I$ are precisely the simple roots of the unipotent radical ${\bf U}_I$
of ${\bf P}_I.$ Thus we have 
$${\bf P}_\Delta ={\bf G}\; \mbox{ and } \;{\bf P}_\emptyset = {\bf P}$$ 
as extreme cases.
Moreover, we have for each subset  $I \subset \Delta$ a unique Levi subgroup ${\bf M}_I$ of ${\bf P}_I$
which contains ${\bf M}.$ 
Let
$$\Phi_I\supset \Phi_I^+ \supset I $$
be its set of roots, positive roots, simple roots with respect to
$\bf{S} \subset \bf{M_I}\cap \bf{P}.$
We denote by $$W=N({\bf S})/Z({\bf S})$$ the relative Weylgroup of ${\bf G}.$ For any subset $I\subset \Delta,$ let  $W_I$ be the parabolic
subgroup of $W$ which is generated by the reflections associated to
$I.$ It coincides with the Weylgroup of ${\bf M}_I.$
Thus we have  
$$W_\Delta=W \mbox{ and } W_\emptyset = \{1\}.$$ 
If ${\bf H}$ is any linear algebraic group defined over $F$, then we denote by $X^\ast({\bf H})_{F}$
its group of $F$-rational characters.

Whereas we denote algebraic groups defined over $F$ by boldface
letters, we use ordinary letters for their groups
$$G:={\bf G}(F),\; P_I:={\bf P}_I(F),\; M_I:={\bf M}_I(F) , \ldots$$ 
of $F$-valued points. We supply these groups with the canonical topology given by $F.$
These are locally profinite topological groups. 
Let ${\bf M \subset G}$ be a Levi subgroup.
Put 
$$^0M= \bigcap_{\alpha \in X^\ast({\bf M})_{F}} kern|\alpha|_{\R}.$$  
This is a normal open subgroup generated by all compact subgroups of $M$ (cf. \cite{BW} ch. X 2.2).
Moreover, the quotient $M/^0M$ is a finitely generated free abelian group of rank equal to the $F$-rank of $Z(\bf{M}).$ 
The valuation map gives rise to a natural homomorphism of groups
\begin{eqnarray}
\Theta_M: X^\ast({\bf M})_{F} \longrightarrow Hom (M/^0M, \Z)
\end{eqnarray}
defined by $\Theta(\chi) = val \circ \chi(F),$ where
$\chi(F):M\rightarrow F^\times$ is the induced homomorphism on
$F$-valued points.  It is easily seen that $\Theta_M$ is
injective. Further the source and the target of $\Theta_M$ are both free
$\Z$-modules of the same rank. Thus we may identify  $X^\ast({\bf M})_{F}$ as a lattice in $Hom (M/^0M, \Z).$ 

We fix a self-injective ring $R,$ i.e., $R$ is an injective object in
the category $Mod_R$ of $R$-modules.  Let $i:\Z \rightarrow
R$ be the canonical homomorphism. Then we have $ker(i)=d\Z,
$ for some integer $d\in \N.$ We suppose that $R$ fulfills the following assumptions.
\begin{enumerate}
\item  The pro-order $|G|$ of $G$ is invertible in $R,$ i.e.,
  $|G|$ is prime to $d$ (see \cite{V1} for the definition of the pro-order).
In particular $i(q)\in R^\times.$  
\item Let 
$$\rho = \det Ad_{Lie({\bf U})}|{\bf S} \in X^\ast({\bf S})_F$$  be the character given by the determinant of the
adjoint representation of ${\bf P}$ on $Lie({\bf U})$ restricted to ${\bf S}.$  Write $\rho$ in the shape
$$\rho = \sum_{\alpha \in \Delta} n_\alpha \alpha,$$ 
where $n_\alpha \in \N.$ Following the definition of  an  algebraically
closed field which is {\it bon} for $G$  (see \cite{D}),
we impose on $R$ that $d$ is prime to 
$$\prod_{r\leq sup\{n_\alpha; \;\alpha \in \Delta\}} (1-q^r).$$
\item Let $E/F$ be a finite Galois splitting
field of ${\bf G}.$ Then we further suppose that $d$ is prime to
the order of the Galois group $Gal(E/F),$ i.e, $i(|Gal(E/F)|) \in R^\times.$
\item Finally we assume that the monomorphism $\Theta_{M_I}$  becomes an isomorphism after base change to $R$
for all $I\subset \Delta.$
\end{enumerate}

\noindent {\bf Remarks:} 
(1) Examples of such rings are given by fields of characteristic zero or
by $R = \Z/n\Z$ with $n\in \N$ suitable chosen.\\
\noindent (2) If $R$ is an algebraically
closed field, then  condition 1 corresponds to the case {\it banal}
in the sense of Vign\'eras (see \cite{V1}).

Suppose for the moment that $G$ is an arbitrary locally profinite
group. We agree that all $G$-representations  (sometimes we use the
term G-module as well) in this paper are defined over $R.$
Recall that a smooth $G$-representation is a representation $V$ of $G$ such that each $v\in V$ is fixed by a compact subgroup $K\subset G.$
We denote the category of smooth representations by $Mod_G.$
If $V$ is a smooth $G$-module, then we let $\widetilde{V}$ be its smooth dual.
Any closed subgroup $H$ of $G$ gives rise to functors  
$$i^G_H,\; c\mbox{-}i^G_H: Mod_H \rightarrow Mod_G$$
called  the (unnormalized) induction resp. induction with compact support. We recall their definitions.
Let $W$ be a smooth $H$-representation. Then we have
\begin{eqnarray*}
i^G_H(W)& := &\Big\{f:G\rightarrow W;\;  f(hg)=h\cdot f(g) \;\forall h\in H, g\in G, \exists \mbox{ compact }\\
& &   \mbox{ open subgroup } K_f \subset G \mbox{ s.t. } f(gk)=f(g) \; \forall g\in G,k\in K_f \Big\}
\end{eqnarray*}
resp.
\begin{eqnarray*}
c\mbox{-}i^G_H(W)& := &\Big\{f\in i^G_H(W); \mbox{the support of $f$ is compact modulo } H \Big\}.
\end{eqnarray*}
Note that we have $$i^G_H=c\mbox{-}i^G_H,$$ if $H\backslash G$ is compact. If furthermore $W$ is admissible, 
i.e., $W^K$ is of finite type over $R$ for all compact open subgroups
$K\subset G,$ then $i^G_H(W)$ is admissible as well (loc.cit., I,
5.6). Finally, we denote for any $G$-module $V$ by $V^G$ resp. $V_G$ the invariants resp. the coinvariants of $V$ with respect to $G.$

Next, we want to recall the definition of the generalized Steinberg
representations. Let ${\bf 1}$ be the trivial representation of any locally profinite group.
For a subset $I\subset \Delta,$
let 
$$i^G_{P_I}:= i^G_{P_I}({\bf 1})=c\mbox{-}i^G_{P_I}({\bf 1})=C^\infty(P_I\backslash G,R)$$ 
be the admissible representation of locally constant functions on
$P_I\backslash G$ with values in $R.$ If $\Delta \supset J\supset I$
is another subset, then there is an injection
$i^G_{P_J} \hookrightarrow i^G_{P_I}$ which is induced by the natural
surjection $P_I\backslash G \rightarrow P_J\backslash G.$ The generalized Steinberg
representation of $G$ with respect to $I\subset \Delta$ is defined to be the quotient
$$v^G_{P_I}:= i^G_{P_I}/\sum_{I \subset J \subset \Delta\atop J\neq I} i^G_{P_J}.$$ 
In the case $R=\C$ it has been shown that the generalized Steinberg
representations are irreducible and not pairwise isomorphic for
different $I\subset \Delta$ (cf. \cite{Ca2} Thm 1.1). This result has
been generalized by J.-F. Dat \cite{D} to the case of an algebraically closed field
which is {\it bon} and {\it banal} for $G.$ 

We finish this section with introducing some more notations.
We fix a normalized left-invariant  $R$-valued Haar measure $\mu$ on $G$
with respect to a maximal compact open subgroup of $G.$
The existence of such a Haar measure is guaranteed 
by assumption (1) on $R$ (see \cite{V1} I, 2.4). Further, we denote by 
$|\;\; |:F \rightarrow R$ the 'norm' given by the composition
of
\begin{eqnarray*}
F & \longrightarrow &  q^\Z  \\
x & \mapsto & q^{-val(x)}
\end{eqnarray*}
together with the natural homomorphism $\Z[\frac{1}{q}] \rightarrow R.$
Finally, if ${\bf H}$ is any linear algebraic group over $F$, then we put 
$$X({\bf H}):=X^\ast({\bf H})_{F} \otimes_\Z R.$$

\section{The computation}
Let $G$ be an arbitrary locally profinite group which satisfies assumption
1 on $R.$  We want to recall that the
category $Mod_G$ of smooth $G$-representations has then enough injectives
and projectives \cite{V1}. 
This fact provides two different choices for the computation of the Ext-groups 
$Ext^\ast_G(V,W),$ for a given pair of smooth $G$-representations $V,W.$ Notice that 
$$H^i(G,V)=Ext^i_G({\bf 1}, V)$$ is the $i^{th}$ right derived functor of
$$Mod_G \rightarrow Mod_R$$ 
$$V\mapsto V^G,$$ 
whereas
$H_i(G,V)$ denotes the $i^{th}$ left derived functor of the right
exact functor 
$$Mod_G \rightarrow Mod_R$$ 
$$V\mapsto V_G.$$ 
Since $R$ is self-injective, it is easy to see that there is an isomorphism
$$H_i(G,V)^{\vee}=Ext^i_G(V,{\bf 1})$$ 
for all smooth $G$-representations $V$ and for all $i\geq 0$. 
Here the symbol $^{\vee}$ indicates the $R$-dual space.

For our proof of Theorem 1, we need some statements on the
cohomology of smooth representations of locally profinite groups with values in $R.$
Up to Lemma 14 all the statements are well-known in the classical case, i.e., where $R=\C.$ Their proofs  in our situation are essentially  the same. 
But for being on the safe side, we are going to reproduce the arguments
shortly. Up to Lemma 7 - except of Lemma 4 - $G$ is an arbitrary
locally profinite group satisfying  assumption 1 on $R.$

\begin{Lemma}
Let $K\subset G$ be an open compact subgroup. Then $i^G_K({\bf 1})$ is
an injective object in $Mod_G$.
\end{Lemma}

\proof By \cite{V1} I,  4.10 we know that the trivial
$K$-representation ${\bf 1}$ is an injective object. Since the induction functor respects injectives 
(loc.cit. I,  5.9 (b)), we obtain the claim.
\qed

Let $Y$ be the Bruhat-Tits building of ${\bf G}$ over $F.$ We denote by $C^q(Y),\; q\in \N,$ the space of $q$-cochains on $Y$ with values in $R.$
As in the classical case we have the following fact:

\begin{Lemma}
The natural chain complex 
$$0 \rightarrow  R \rightarrow C^0(Y) \rightarrow C^1(Y) \rightarrow \dots \rightarrow  C^q(Y) \rightarrow \dots.$$
is an injective resolution of the trivial $G$-representation ${\bf 1}$  by smooth $G$-modules.
\end{Lemma}

\proof The proof coincides with the proof of \cite{BW} ch. X 1.11  which uses Lemma 3 and the contractibility of the
Bruhat-Tits building  $Y.$  \qed

Our next lemma deals with the Hochschild-Serre spectral sequence. Let $N\subset G$ be a closed subgroup. As it has been pointed out by Casselman
in \cite{Ca2}, the restriction functor from the category of smooth $G$-modules to that of $N$-modules does not preserve injective objects.
For this reason, the standard arguments for proving the existence of the Hochschild-Serre spectral sequence - as in the case of cohomology theory of groups - breaks down. Nevertheless, the restriction functor
preserves projective objects giving a homological variant of the Hochschild-Serre spectral sequence (see appendix of \cite{Ca2}). 

\begin{Lemma}
Let $N \subset G $ be a closed normal subgroup of $G.$ If $V$ is a projective $G$-module, then $V_N$ is a projective $G/N-$module. Thus we 
get for every pair of smooth $G$-modules $V,W,$ such that $N$ acts trivially on $W,$  a spectral sequence
$$E_2^{p,q}=Ext^q_{G/N}(H_p(N,V),W) \Rightarrow Ext^{p+q}_{G}(V,W).$$
If furthermore $N$ resp. $G/N$ is compact, then we have 
$$ Ext^q_{G/N}(V_N,W)= Ext^q_G(V,W)\; \forall q\in \N,$$
resp.
$$Ext^0_{G/N}(H_p(N,V),W) = Ext^p_G(V,W) \; \forall p\in \N .$$
\end{Lemma}

\proof The proof is the same as in the classical case \cite{Ca2} A.9. It starts with the observation that the coinvariant functor is 
left adjoint to the trivial (exact) functor viewing a smooth
$G/N$-module as  a smooth $G$-module. Therefore, $V_N$ is a projective
$G/N$-module, if $V$ is projective.
By \cite{V1} I,  5.10 we know that the restriction functor preserves
projectives. Using the standard-arguments applied to the Grothendieck
spectral sequence, we obtain the first part of the claim.
The reason for the second part is the exactness of the coinvariant resp. fixed-point functor for a 
compact subgroup \cite{V1} I,  4.6. \qed

\begin{Lemma}  Let $V$ and $W$ be smooth representations of $G.$ Suppose that
$W$ is admissible. Then there are  isomorphisms
$$Ext^i_G(V,W)\cong Ext^i_G(\widetilde{W},\widetilde{V}),\; \forall i\geq 0 .$$
\end{Lemma}

\proof Let $$0 \leftarrow V \leftarrow P^0 \leftarrow P^1 \leftarrow  \cdots$$ be a projective resolution of $V.$ Since $R$ is 
self-injective, we conclude as in \cite{V1} I, 4.18 that the functor
$W\mapsto  \widetilde{W}$ from the category of smooth
$G$-representations to itself is exact. 
By \cite {V1} I, 4.13 (2)  we see that the modules $\widetilde{P}^j,
j\geq 0,$ are injective objects in $Mod_G.$
Hence, we obtain an injective resolution  
$$0\rightarrow \widetilde{V} \rightarrow \widetilde{P}^0 \rightarrow \widetilde{P}^1 \rightarrow \dots $$ of $\widetilde{V}.$ 
Moreover, we know by \cite{V1} I,  4.13 (1) that 
$$Hom_G(V,\widetilde{W}) = Hom_G(W,\widetilde{V}),$$ for any pair of smooth $G$-modules
$V,W.$ Since $W$ is admissible, we have 
$W=\widetilde{\widetilde{W}}$ (see \cite{V1} 4.18 (iii)) and the claim follows. \qed

\noindent In the special case $W={\bf 1}$ we obtain:

\begin{Corollary}
Let $V$ be a smooth representation of $G.$  Then there are isomorphisms
$$H^i(G,\widetilde{V})\cong H_i(G,V)^\vee, \; \forall i\geq 0.$$
\end{Corollary}

From now on, we suppose again that $G$ is the set of $F$-valued points
of some reductive algebraic group defined over $F.$

\begin{Lemma}  Let $Q\subset G$ be a parabolic subgroup with Levi
  decomposition $Q=M\cdot N.$
Let $V$ resp. $W$ be a smooth representation of $G$ resp. $M.$ Extend $W$ trivially to a representation of $Q.$
Then we have for all $i\geq 0$ isomorphisms
$$Ext^i_G(V,i^G_Q(W))\cong Ext^i_M(V_N,W) .$$
\end{Lemma}

\proof By Frobenius reciprocity \cite{V1} I, 5.10 we deduce that 
$$Ext^\ast_G(V,i^G_Q(W)) = Ext^\ast_Q(V,W).$$ Since $N$ is a union of
open compact subgroups, we deduce from \cite{V1} I, 4.10 
the exactness of the functor
$$Mod_G \rightarrow Mod_R$$ 
$$W\mapsto W_N.$$  
Thus the statement follows from Lemma 5. \qed

After having established the main techniques for computing 
cohomology of representations, we are able to take the first step
in order to proof Theorem 1. The following proposition is also well-known in the classical case.

\begin{Proposition}
We have 
$$H^\ast(G,{\bf 1})= \Lambda^\ast X({\bf G }),$$ 
where $\Lambda^\ast X({\bf G})$ denotes the exterior algebra of
$X({\bf G})$.
\end{Proposition}

\proof We copy the proof of the classical case \cite{BW} Prop. 2.6, ch. X.

$1^{st}$ case: ${\bf G}$ is semi-simple and simply connected. Then we apply the $G$-fixed point functor to the resolution of the trivial
representation in Lemma 4. The result is a constant coefficient system
on a base chamber inside the Bruhat-Tits building, which is contractible. Thus, we obtain $H^\ast(G,{\bf 1})= H^0(G,{\bf 1})= R.$

$2^{nd}$ case: ${\bf G}$ is semi-simple. Then we consider its simply
connected covering ${\bf G'}\rightarrow {\bf G}.$ The induced
homomorphism $G' \rightarrow G $ has finite kernel, its image is a
closed cocompact normal subgroup. We apply Lemma 5 to
$G',$ $\sigma(G')$ and $ N:=ker(G' \rightarrow G).$

$3^{rd}$ case: ${\bf G}$ is arbitrary reductive. 
Let $D{\bf G}$ be the derived group of ${\bf G}$ and put $G'=D{\bf G}(F).$
Then we  have $G\supset{} ^0G \supset DG'.$ Moreover, the quotient $^0G / DG'$ is compact,
where $DG'$ denotes the derived group of $G'.$ Therefore, we conclude by
the previous case, Lemma 5 and Corollary 7 that 
$$H^\ast(^0G,{\bf 1})= H^\ast(DG',{\bf 1}) =H^0(DG',{\bf 1})=R.$$
With the same arguments, we see that
$$H^\ast(G,{\bf 1})= H^\ast(G/^0G,{\bf 1}).$$
Now it is known that the cohomology of a finite rank free commutative (discrete) group $L$  coincides with the cohomology of the corresponding
torus: 
$$H^\ast(L,{\bf 1})= \Lambda^\ast(Hom(L,\Z))\otimes_{\Z} R.$$  Applying this fact to $G/^0G$, we  get
$$H^\ast(G,{\bf 1}) = \Lambda^\ast(Hom(G/^0G,\Z))\otimes_{\Z} R.$$
By assumption 4 on $R$ we have  $Hom(G/^0G,\Z)\otimes_{\Z} R \cong X({\bf G})$
from which the result follows. \qed

\begin{Corollary}
Let $I\subset \Delta.$ Then we have
$$H^\ast(G,i^G_{P_I})= H^\ast(P_I,{\bf 1})=H^\ast(M_I,{\bf
  1})=\Lambda^\ast X({\bf M}_I).$$
\end{Corollary}

\proof The statement follows from Lemma 8, Proposition 9 and by our assumption 4  on $R$.
\qed

In order to compute the cohomology of generalized Steinberg
representations, we need the following proposition.
For two subsets $I\subset I' \subset \Delta$ with $|I'\setminus I| =1$, we
let 
$$p_{I,I'}:i^G_{P_{I'}} \longrightarrow i^G_{P_I}$$
be the natural homomorphism induced by the surjection $G/P_I \rightarrow G/P_{I'}.$
For arbitrary subsets  $I,I' \subset \Delta$, with $|I'| -|I|=1$ 
and $I'=\{\beta_1,\ldots,\beta_r\},$
we put $$d_{I,I'}=\left\{ \matrix{ (-1)^i p_{I,I'} & I' =
I  \cup \{\beta_i\} \cr 0 & I \not\subset I' } \right. .$$

\begin{Proposition}
Let $I\subset \Delta.$ The complex
$$0 \rightarrow  i^G_G \rightarrow \bigoplus_{I\subset K \subset \Delta \atop |\Delta\setminus K|=1}i^G_{P_K} \rightarrow \bigoplus_{I\subset K \subset \Delta \atop |\Delta\setminus K|=2}i^G_{P_K} 
\rightarrow \dots \rightarrow \bigoplus_{I\subset K \subset \Delta
  \atop |K\setminus I|=1}i^G_{P_K}\rightarrow i^G_{P_I} \rightarrow
v^G_{P_I}\rightarrow 0,$$
with differentials induced by the $d_{J,J'}$ above is acyclic.
\end{Proposition}

\proof See  Prop. 13, \S 6 of \cite{SS} for the case of
 $I=\{\alpha_1,\alpha_2,\ldots,\alpha_i\},\,i\geq 1,$ and $G=GL_n.$
 The proof there is only formulated for coefficients in the ring  of integers $\Z.$ However, the 
proof holds for arbitrary rings, since it is of combinatorial nature. 
\\A different approach consists of using Proposition 6 
of \S 2 in \cite{SS}. It says: Let $G_1,\ldots,G_m$ be a family of subgroups in some bigger group $G.$
Suppose that the following identities are satisfied for all subsets $A,B \subset \{1,\ldots, m\}$:
$$(\sum_{i\in A} G_i) \cap (\bigcap_{j\in B} G_j) = \sum_{i \in A}(G_i\cap (\bigcap_{j\in B} G_j)).$$
Then the natural (oriented) complex
$$G \leftarrow \bigoplus_{i=1}^m G_i \leftarrow \bigoplus_{i,j=1 \atop i<j}^ m  G_i \cap G_j \leftarrow \bigoplus_{i,j,k=1\atop i<j<k}^m G_i \cap G_j \cap G_k\leftarrow \cdots$$
is an acyclic resolution of $\sum_i G_i \subset G.$
We apply this proposition to the $G$-modules $i^G_{P_K},$ where $I\subset K
\subset \Delta$ and $|\Delta \setminus K|=1.$ The condition of the
proposition is fulfilled. Indeed,
we have
$$i^G_{P_I}\; \cap\; i^G_{P_J} = i^G_{P_{I \cup J}}$$  
and
$$ i^G_{P_I}\; \cap\; (i^G_{P_J} + i^G_{P_K})  = (i^G_{P_I}  \cap
i^G_{P_J})\; +\;  (i^G_{P_I}  \cap i^G_{P_K}), $$
for all subsets $I,J,K \subset \Delta.$
The first identity follows from the fact that $P_{I\cup J}$ is the parabolic subgroup generated by $P_I$ and $P_J.$ 
For the second one confer \cite{BW} 4.5, 4.6  resp. \cite{L} 8.1,
8.1.4 (The statement there is formulated in the case where $R = \C.$ 
The result holds also in our general situation. The proof relies on the exactness
of the Jacquet-functor and a description of the $S$-modules
$(i^G_{P_I})_U$ using the filtration in the proof of
Proposition 15).
\qed

\begin{Theorem} Let $G$ be semi-simple and let $I\subset \Delta.$ Then we have
$$H^i(G,v^G_{P_I})=\left\{ \begin{array}{r@{\quad:\quad}l} R & i =|\Delta \setminus I| \\ 0 & \mbox{otherwise} \end{array} \right.$$
\end{Theorem}

\proof The proof is the same as in  Prop. 4.7, ch. X of \cite{BW}. A not very different approach works as follows. 
Apply the cohomology functor $H^\ast(G,-)$ to the acyclic complex of
Proposition 11. We obtain a complex
$$0 \rightarrow \Lambda^\ast X({\bf G}) \rightarrow 
\bigoplus_{I\subset K \subset \Delta \atop |\Delta\setminus K|=1}
\Lambda^\ast X({\bf M}_K) \rightarrow  \dots  \rightarrow \bigoplus_{I\subset K \subset \Delta \atop
  |K\setminus I|=1} \Lambda^\ast X({\bf M}_K)
\rightarrow  \Lambda^\ast X({\bf M }_I) \rightarrow 0 .$$ 
Using the Hochschild-Serre spectral sequence, we may assume without
loss of generality that  ${\bf G}$ is simply connected. Suppose that
${\bf G}$ is split. In this case it is well-known (cf. \cite{J} ch II,
1.18 ) that $X^\ast({\bf M}_K)_{F}$ may be identified with the submodule
of $X^\ast({\bf S})_{F}$
defined by $$\{ \chi \in X^\ast({\bf S})_{F};\; \langle \chi,
\alpha^\vee \rangle\; =0 \;\forall \alpha \in K\},$$
where $\langle\,\; ,\; \rangle: X^\ast({\bf S})_{F} \times X_\ast({\bf
  S})_{F} \rightarrow \Z$ is the natural pairing.
If we denote by $\{\omega_\alpha \in X^\ast({\bf S})_{F};\; \alpha \in
\Delta\}$ the fundamental weights of {\bf G} with respect to 
${\bf S} \subset {\bf  P }$, then we get
$$X({\bf M}_K) \cong \bigoplus_{\alpha\in  \Delta \setminus K} R\cdot\omega_\alpha \subset X({\bf S}).$$ 
Thus we see - again by using Prop. 6, \S 2 of \cite{SS} -   that the
complex above is acyclic with respect to $\Lambda^r$ for 
$$r< rk(Z({\bf M}_I)) = |\Delta\setminus I|.$$ 
In the case $rk(Z({\bf M}_I))=r$ all the entries of the complex vanish
except of $\Lambda^r X({\bf M}_I)=R.$
\\In the general case, let $E/F$ be our fixed Galois splitting field of ${\bf G}.$ Then we deduce with the same arguments that the corresponding 
complex of $E$-rational characters has the desired property. Applying
the $Gal(E/F)$-fixed point functor to this complex  
yields the claim. Note that the fixed point functor is exact
by assumption 3 on $R.$\qed

\noindent For attacking Theorem 1 we still need two lemmas.

\begin{Lemma} Let  $V$ be a smooth representation of $G.$ Suppose that
  there exists an element $z \in Z(G)$ in the center
of $G$ and an element $c \in R,$ such that $ c-1 \in R^\times $ and
$z\cdot v =c \cdot v$ for all $v\in V.$  Then we have
$$H^\ast(G,V)=0 .$$
\end{Lemma}

\proof See Prop. 4.2, ch. X \cite{BW} for the classical case. We repeat shortly the argument. By identifying 
Ext-groups with the Yoneda-Ext-groups, we have to show that for all $n
\in \N$, all $n$-extensions of ${\bf 1}$ by $V$ are trivial. 
More generally, we will show that
if $U$ is a $R$-module with trivial $G$-action, then there are no non-trivial extensions of $U$ by $V.$ In fact, let
$$E^\bullet: 0 \rightarrow V \rightarrow E^1 \rightarrow E^2 \rightarrow \cdots \rightarrow E^n \rightarrow U \rightarrow 0$$ be an 
arbitrary $n-$extension. 
Since $z$ lies in the center of $G,$
it defines an endomorphism of $E^\bullet$ and we get the identity
$E^\bullet = c.E^\bullet.$ Here $c.E^\bullet$ denotes the scalar multiplication
of $R$ on the module $Ext^n_G(U,V)$ (confer \cite{M} ch. III, Theorem 2.1). Thus, we have $0=E^\bullet - c.E^\bullet = (1-c).E^\bullet.$ Since 
$1-c \in R^\times,$ we conclude that $E^\bullet=0 \in  Ext^n_G(U,V).$ \qed

\begin{Lemma} 
Let $H\subset G$ be a closed subgroup and  let $W$ be a smooth representation of $H.$ Then we have
$$ \widetilde{c-i^G_H(W)} \cong i^G_H(\widetilde{W} \delta_H) ,$$
where $\delta_H$ is the modulus character of $H.$
\end{Lemma}

\proof  
This follows from \cite{V1} I,  5.11 together with the fact that $G$ is unimodular. \qed

\begin{Proposition} Let $G$ be semi-simple and let $I,J \subset \Delta.$ Then we have
$$Ext^\ast_G(i^G_{P_I},i^G_{P_J}) =\left\{
  \begin{array}{r@{\quad:\quad}l}\Lambda^\ast X({\bf M}_J)  & \mbox{ if } J\subset I \\ 0 & \mbox{ otherwise } \end{array} \right.$$
\end{Proposition}

\proof By Lemma 8  we have for all $i\geq 0$ isomorphisms
$$Ext^i_G(i^G_{P_I},i^G_{P_J})\cong Ext^i_{M_J}((i^G_{P_I})_{U_J},{\bf 1}),$$ 
where $(i^G_{P_I})_{U_J}$ is the Jacquet-module of $i^G_{P_I}$ with respect to $M_J.$
In the case $R=\C$  there is constructed in \cite{Ca3} 6.3 -  a
substitute for the Mackey formula - a decreasing $\N$-filtration
$\F^\bullet$ of smooth $P_J$-submodules on $i^G_{P_I}$  defined by
$$\F^i=\{f\in i^G_{P_I}\;;\, supp(f) \subset \bigcup_{w\in W_I\backslash
  W/W_J \atop l(w)\geq i} P_I\backslash P_IwP_J\},\; i\in \N.$$
Here the length $l(w)$ of a double coset  $w\in W_I\backslash W /W_J$ is  the length of its
  Kostant-representative which is the one of minimal length
  within its double coset. In the following we will identify the
  double cosets with its Kostant-representatives.
There are a canonical isomorphisms 
$$gr^i_{\F^\bullet}(i^G_{P_I})\cong \bigoplus_{w \in W_I\backslash W/ W_J \atop l(w)=i} c-i^{P_J}_{P_J \cap w^{-1}P_I w},$$
for all $i\geq 0.$
Furthermore, we have for every $w\in W_I\backslash W/ W_J $ an isomorphism 
$$ (c\mbox{-}i^{P_J}_{P_J \cap w^{-1}P_I w})_{U_J} \cong c\mbox{-}i^{M_J}_{M_J \cap w^{-1}P_Iw}(\gamma_w),$$
where $\gamma_w$ is the modulus character of $P_J \cap w^{-1}P_I w$ acting on $U_J/U_J\cap w^{-1}P_Iw.$ 
The first isomorphism is a corollary of Prop. 6.3.1 (loc.cit.) (see
also  \cite{V1} I, 1.7 (iii)), 
whereas the second one is the content of Prop. 6.3.3 (loc.cit.).
In the general case, i.e., for our specified ring $R$, the same formulas hold,
since the proof can be taken over word by word. Since $M_J \cap w^{-1}P_Iw$ is 
a parabolic subgroup in $M_J,$
we observe that $c$-$i^{M_J}_{M_J \cap w^{-1}P_Iw}(\gamma_w)=i^{M_J}_{M_J \cap w^{-1}P_Iw}(\gamma_w).$ 
From the definition we see that $\gamma_w$ is the norm of the rational character
$$\det Ad_{Lie({\bf U_J})}/\det Ad_{w^{-1}Lie({\bf P_I})w\cap Lie(\bf{U_J})} \in
X^\ast({\bf P_J} \cap w^{-1}{\bf P_I}w) .$$ Its restriction to $S$ is given by
\begin{eqnarray}
{\gamma_w}_{ |S} = |\prod_{\alpha \in \Phi^+\setminus \Phi^+_J \atop w\alpha \in \Phi^-\setminus \Phi^-_I} \alpha |.
\end{eqnarray}
Fix an element $w \in W_I\backslash W/ W_J.$ We are going to show that 
$$Ext^\ast_{M_J}(i^{M_J}_{M_J \cap w^{-1}P_I w}(\gamma_w),{\bf 1})=0,$$
unless $w=1$ and $J\subset I.$ Since the Jacquet-functor is exact, this will give by successive application of the long exact cohomology 
sequence with respect to the filtration $\F^\bullet$ the statement of our proposition.
By Lemma 6 and Lemma 14 we conclude  that
$$Ext^\ast_{M_J}(c\mbox{-}i^{M_J}_{M_J \cap w^{-1}P_I w}(\gamma_w),{\bf 1}) 
\cong Ext^\ast_{M_J}({\bf 1},i^{M_J}_{M_J \cap w^{-1}P_I w}(\tilde{\gamma}_w\delta_{M_J \cap w^{-1}P_I w})),$$ where
$\delta_{M_J \cap w^{-1}P_I w}$ is the modulus character of the parabolic subgroup \mbox{$M_J \cap w^{-1}P_I w$}
of $M_J$ and $\tilde{\gamma}_w$ is the smooth dual of $\gamma_w.$  The Levi decomposition of the latter group is given by
$$M_J \cap w^{-1}P_I w=M_{J \cap w^{-1}I}\cdot (M_J \cap w^{-1}U_I w)$$
(see \cite{C} Prop. 2.8.9.).
So, the restriction of $\delta_{M_J \cap w^{-1}P_I w}$ to $S$ is
the norm of the rational character 
$$\prod_{\alpha \in \Phi^+_J \atop w\alpha \in \Phi^+\setminus
  \Phi^+_I}\alpha,$$
i.e.,
\begin{eqnarray} 
{\delta_{M_J \cap w^{-1}P_I w}}_{|S}= |\prod_{\alpha \in \Phi^+_J \atop w\alpha \in \Phi^+\setminus \Phi^+_I}\alpha|.
\end{eqnarray}
In the case where $J\not\subset I$ or $w\neq1$ we deduce from the following lemma 
the existence of an element $z$ in the center of $M_{J\cap w^{-1}I},$
such that
$$\tilde{\gamma}_w(z)\delta_{M_J \cap w^{-1}P_I w}(z)-1 \in R^\times.$$
By Lemma 13 we conclude that   
$$Ext^\ast_{M_J}(c\mbox{-}i^{M_J}_{M_J \cap w^{-1}P_I w}(\gamma_w),{\bf 1})=0.$$   
In the case $J\subset I$ we obtain therefore an isomorphism 
$$Ext^\ast_G(i^G_{P_I},i^G_{P_J}) \cong Ext^\ast_{M_J}({\bf 1},{\bf
  1})=\Lambda^\ast X({\bf M}_J)$$ which is induced by the element $w=1.$
\qed

\begin{Lemma} Let $J\not\subset I$ or $w\neq1$. Then
there exists an element $z\in Z(M_{J\cap w^{-1}I})$ such that $\tilde{\gamma}_w(z) \delta_{M_J \cap w^{-1}P_I w}(z) - 1 \in R^\times.$
\end{Lemma}

\proof 
$1^{st}$ case: Let $w \neq 1.$ Then we have $\gamma_w \neq {\bf 1}.$ In fact,
$\gamma_w=1$ would imply that  
$$Lie(U_J) \subset Lie(w^{-1}P_Iw)$$ or equivalently  $U_J\subset w^{-1}P_Iw.$
But in general one has $$P_{J\cap w^{-1}I}= (P_J\cap w^{-1}P_Iw)\cdot U_J$$ (\cite{C}, Prop. 2.8.4).
Thus, we deduce that the intersection  $P_J\cap w^{-1}P_Iw$ is a
parabolic subgroup. This is only true if $w=1.$
\\We want to recall that for any subset $K\subset \Delta$
the maximal split torus in the center $Z(M_K)$ of $M_K$ coincides with
the connected component of the identity in
$\bigcap_{\alpha \in K}  kern(\alpha) \subset S.$ 
Since the center of $M_J$ is contained in $M_{J\cap w^{-1}I},$ it is  
enough to construct an element  $z\in Z(M_J)$ which has the desired property.
From the representation (2) we may easily conclude the existence of  an element  $z\in Z(M_J)$
with $\tilde{\gamma}_w(z)\neq 1.$ Our purpose is to show the
existence of an element $z\in Z(M_J)$ such that
$\tilde{\gamma}_w(z) -1 \in R^\times.$
We may suppose that ${\bf G}$ is adjoint. Let
$$\{\omega_\alpha \in X_\ast(S);\; \alpha \in \Delta\}$$ be the dual
base (co-fundamental weights) of $\Delta,$ i.e.,
$\langle \omega_\beta , \alpha \rangle = \delta_{\alpha,\beta},$
for all $\alpha,\beta\in \Delta.$ Since $\gamma_w \neq {\bf 1}$ it is
possible to
find a root $\alpha \in \Delta\setminus J$  such that $w\alpha \in
\Phi^-\setminus \Phi^-_I.$ 
Put $$z:=\omega_\alpha(\varpi_F^{-1}).$$ Then we have $z\in
Z(M_J)$ and  
$$\tilde{\gamma}_w(z) -1 = q^r - 1 $$ for some $1\leq r \leq n_\alpha.$
By assumption 2 on $R$ we know that the product $\prod_{r\leq sup\{n_\alpha; \;
\alpha \in \Delta\}} (1-q^r)$ is invertible in $ R.$
Further we see from the expression (3) that  $\delta_{M_J\cap w^{-1}P_I w}(z)=1.$ This
gives the proof in the first case. 

\noindent $2^{nd}$ case: Let  $w= 1$ and $J\not\subset I.$ Then we have
$\gamma_w=1.$ Since $J\not\subset I$,  we see that the restriction of $\delta_{M_J\cap P_I}$ to  $Z(M_{J\cap I})$ is not trivial. 
Again, we can find as in the first case an
element $z \in Z(M_{J \cap I})$ such that $\delta_{M_J\cap P_I}(z) - 1  \in R^\times.$
\qed

\begin{Proposition} Let $G$ be semi-simple and let $I,J \subset \Delta.$ Then we have
$$Ext^\ast_G(v^G_{P_I},i^G_{P_J})=\left\{
  \begin{array}{r@{\quad:\quad}l} \Lambda^\ast X({\bf M}_J)[-|\Delta\setminus I|]  & \Delta=I\cup J \\ 0 & \mbox{otherwise} \end{array} \right.$$
\end{Proposition}

\proof We apply the acyclic complex of Proposition 11 to the representation $v^G_{P_I}.$ This yields a double complex
$$ 0\rightarrow Ext_G^\ast(i^G_{P_I},i^G_{P_J}) \rightarrow \bigoplus_{I \subset L \subset \Delta \atop |L\setminus I|=1} Ext_G^\ast(i^G_{P_L},i^G_{P_J})
\rightarrow \bigoplus_{I \subset L \subset \Delta \atop |L\setminus I|=2} Ext_G^\ast(i^G_{P_L},i^G_{P_J})
\rightarrow \dots $$
$$ \dots \rightarrow \bigoplus_{I \subset L \subset \Delta \atop |\Delta\setminus L|=1} Ext_G^\ast(i^G_{P_L},i^G_{P_J}) \rightarrow 
Ext_G^\ast(i^G_G,i^G_{P_J})\rightarrow 0,$$
such that its associated spectral sequence converges to $Ext^\ast_G(v^G_{P_I},i^G_{P_J}).$
By Proposition 15 we see  that $K:=I\cup J$ is the minimal subset of $\Delta$ 
containing $I$ with $Ext^\ast_G(i^G_{P_K},i^G_{P_J})\neq 0.$ 
Hence, the double complex reduces to the double-complex
$$0 \rightarrow   \Lambda^\ast X({\bf M}_J) \rightarrow \bigoplus_{K \subset L \subset \Delta \atop |L\setminus K|=1} 
\Lambda^\ast X({\bf M}_J) \rightarrow \bigoplus_{K \subset L \subset \Delta \atop |L\setminus K|=2} 
\Lambda^\ast X({\bf M}_J) \rightarrow \dots $$
$$ \dots \rightarrow \bigoplus_{K\subset L \subset \Delta \atop
  |\Delta\setminus L|=1} \Lambda^\ast X({\bf M}_J) 
\rightarrow  \Lambda^\ast X({\bf M}_J) \rightarrow 0.$$
In the case of $K=\Delta$ we are obviously done. In the case $K\neq \Delta$ we see that the cohomology of the double complex 
vanishes, since it is  a constant coefficient system on the standard simplex corresponding to the set $K.$ 
 \qed

\noindent {\bf Proof of Theorem 1:} This time we apply Proposition 11 to $v^G_{P_J}.$ 
This yields  a double complex
$$ 0\rightarrow Ext_G^\ast(v^G_{P_I},i^G_G) \rightarrow \bigoplus_{J \subset L \subset \Delta \atop |\Delta\setminus L|=1} Ext_G^\ast(v^G_{P_I},i^G_{P_L})
\rightarrow \bigoplus_{J \subset L \subset \Delta \atop |\Delta\setminus L|=2} Ext_G^\ast(v^G_{P_I},i^G_{P_L})
\rightarrow \dots $$
$$ \dots \rightarrow \bigoplus_{J \subset L \subset \Delta \atop |L\setminus J|=1} Ext_G^\ast(v^G_{P_I},i^G_{P_L}) \rightarrow 
Ext_G^\ast(v^G_{P_I},i^G_{P_J})\rightarrow 0,$$
such that its associated spectral sequence converges to $Ext^\ast_G(v^G_{P_I},v^G_{P_J}).$
By Proposition 17 we conclude that the minimal
subset $K$ of $\Delta$ containing $J$ and such that $Ext^\ast_G(v^G_{P_I},i^G_{P_K})\neq 0$ is 
$$K=(\Delta \setminus I) \cup J= (\Delta\setminus I)\; \dot\cup\; (I\cap J).$$
Thus the complex above reduces to
$$0 \rightarrow \Lambda^\ast X({\bf G})[-|\Delta\setminus I|] \rightarrow 
\bigoplus_{K\subset L \subset \Delta \atop |\Delta\setminus L|=1}
\Lambda^\ast X({\bf M}_L)[-|\Delta\setminus I|] \rightarrow \dots $$
$$ \dots \rightarrow \bigoplus_{K\subset L \subset \Delta \atop
  |L\setminus K|=1} \Lambda^\ast X({\bf M}_L)[-|\Delta\setminus I|] 
\rightarrow  \Lambda^\ast X({\bf M}_K)[-|\Delta\setminus I|] \rightarrow 0 .$$ 
This double-complex is precisely - up to shifts -  the double-complex for the computation of the cohomology of $v^G_{P_K},$ 
for a semi-simple group $G$  (cf. Theorem 12 resp. \cite{BW} ch. X, Prop. 4.7) !
Thus, we obtain an isomorphism
$$H^\ast(G,v^G_{P_K})[-(|J| - |K|) - |\Delta\setminus I|]\cong Ext_G^\ast(v^G_{P_I},v^G_{P_J}).$$
It remains to compute the degree $d,$ where the latter space does not vanish.
The degree is by Theorem 12 equal to
\begin{eqnarray*}
d& = & |\Delta\setminus K| + |\Delta\setminus I| + |J| - |K|\\
& = & |\Delta\setminus(\Delta\setminus I \;\dot\cup\; (I\cap J))| + |\Delta\setminus I| + |J| - |\Delta\setminus I\; \dot\cup\; (I\cap J))|  \\
& = & |I\cap \Delta\setminus(I\cap J))| + |J| - |I\cap J| \\
& = & |I\setminus(I\cap J))| + |J| - |I\cap J| \\
& = & |I|  - |I \cap J| + |J| - |I\cap J| \\
& = & |I\cup J| - |I\cap J|. \;\;\;\;\;\;\;\;\;\;\;\;\;\;\;\;\;\;\;\;\;\;\;\;\;\;\;\;\;\;\;\;\;\;\;\;\;\;\;\;\;\;\;\;\;\;\;\;\;\;\;\;\;\;\;\;\;\;\;\;\;\;\;\;\;\; \square
\end{eqnarray*}

\noindent {\bf Remark:} An argument of J.-F. Dat shows
that Theorem 1 even holds if $R$ is not self-injective. In fact, in
his paper \cite{D} Theorem 4.4 he first shows the statement for an
algebraically closed field which is {\it bon} and {\it banal} for $G$.
Then he uses this result to deduce the general case by elementary
commutative algebra.

\noindent {\bf Proof of Corollary 2 :} Consider the projection $G\rightarrow G/Z(G)$ onto the adjoint group of $G.$ 
The action of $Z(G)$ on $v^G_{P_I}$ and $v^G_{P_J}$ is trivial. 
By applying Lemma 5 to this situation we get a spectral sequence
$$Ext_{G/Z(G)}^\ast(H_\ast(Z(G),v^G_{P_I}),v^G_{P_J}) \Rightarrow  Ext_G^\ast(v^G_{P_I},v^G_{P_J}).$$
By the proof of Proposition 9  we deduce that
$$H^\ast(Z(G),{\bf 1})=\Lambda^\ast Hom(Z(G)/^0 Z(G),\Z)\otimes R  \cong \Lambda^\ast R^d.$$
Therefore, we get
$$H^\ast(Z(G),v^G_{P_I})=H^\ast(Z(G),{\bf 1}) \otimes v^G_{P_I} \cong  \bigoplus_{j=0}^d \;\;(v^G_{P_I})^{d \choose j}  .$$
Now we apply Theorem 1 together with Corollary 7. \qed

\bigskip
In the remainder of this paper we give another corollary in
the case of the general linear group and $R=\C.$  This corollary has been
pointed out to me by C. Kaiser.

Let $G=GL_n$ with $n=r\cdot k$ for some integers $k,r>0.$ 
Let $P_{r,k}$ be the upper block parabolic subgroup containing the
Levi subgroup 
$$\underbrace{GL_r \times \cdots \times GL_r}_k.$$
Let $\sigma$ be an
irreducible cuspidal representation of $GL_r.$ For any integer $i\geq 0$ 
we put $\sigma(i)=\sigma\otimes |\det|^i,$ where
$\det:GL_r\rightarrow F^\times$ is the determinant. Consider the graph
$\Gamma$ consisting of the vertices
$\{\sigma,\sigma(1),\ldots,\sigma(k-1)\}$ and the edges
$\{\{\sigma(i),\sigma(i+1)\}; i=0,\ldots,k-2\}.$ Thus we can
illustrate $\Gamma$ in the shape
$$ \sigma - \sigma(1) - \cdots - \sigma(k-1).$$
An orientation of $\Gamma$ is given by choosing a direction on each
edge. Denote by $Or(\Gamma)$ the set of orientations on $\Gamma.$ 

Let ${\cal J}$ be the set of irreducible subquotients of
$\tilde{i}^G_{P_{r,k}}(\sigma\otimes \sigma(1) \otimes \cdots \otimes
\sigma(k-1)),$ 
where $\tilde{i}^G_{P_{r,k}}$ denotes the normalized induction functor.
Following \cite{Z} 2.2, there is a bijection 
$$\omega: Or(\Gamma) \rightarrow {\cal J},$$
which we briefly describe.
Let $S_k$ be the symmetric group in the set \mbox{$\{0,\ldots,k-1\}.$}
Consider the map
$$S_k\rightarrow Or(\Gamma)$$
$$w \mapsto \Gamma(w)$$
defined as follows. The edge $\{\sigma(i),\sigma(i+1)\}$ is oriented from
$\sigma(i)$ to $\sigma(i+1)$ - symbolized as $\sigma(i) \rightarrow
\sigma(i+1)$ - if and only if $w(i) <  w(i+1).$  
On easily verifies the surjectivity of this map.
Let $\vec{\Gamma}$ be an orientation of $\Gamma.$ Choose an element
$w\in S_k$ such that $\vec{\Gamma}= \Gamma(w).$ Then
$\omega(\vec{\Gamma})$
is defined to be the unique irreducible quotient of
$$\tilde{i}^G_{P_{r,k}}(\sigma(w(0))\otimes \cdots \otimes \sigma(w(k-1))).$$
In loc.cit. 2.7  it is shown that this representation does not depend on the chosen representative
$w.$

Denote by $\Delta_k=\{\alpha_0,\ldots,\alpha_{k-2}\}$ the set of simple roots of $GL_k$  with respect
to the standard root system of $GL_k.$ Let ${\cal P}(\Delta_k)$ be its
power set. For a subset $I\subset \Delta_k,$ we let
$\Theta(I) \in Or(\Gamma)$ be the orientation of $\Gamma$ defined by $\sigma(i)
\rightarrow \sigma(i+1)$ if and only if $\alpha_i \in I,
i=0,\ldots,k-2.$  It is easily seen that we get in this way a bijection 
$$\Theta: {\cal P}(\Delta_k) \rightarrow Or(\Gamma).$$
For any subset $I\subset \Delta_k,$ we put 
$$v^G_I(\sigma):=\omega(\Theta(I)).$$ 

\noindent {\bf Example:} Consider the special case $r=1$ and $\sigma=|\;|^{\frac{1-n}{2}}.$
Then we have $P_{r,k}=P,$
$$\tilde{i}^G_P(\sigma \otimes \cdots \otimes
\sigma(n-1))= i^G_P$$ and   
$$v^G_I(\sigma) = v^G_{P_I},$$
for all $I\subset \Delta=\Delta_k.$

\begin{Corollary} Let $I,J \subset \Delta_k.$ Set $i:= |I\cup J| -
  |I\cap J|.$ Then we have
$$Ext^\ast_G(v^G_I(\sigma),v^G_J(\sigma))=  R[-i] \oplus R[-i-1]. $$
\end{Corollary}

\proof We make use of the theory of types of Bushnell and Kutzko
\cite{BK} (see also \cite{V2}). Let $(K,\lambda)$
be the type of the block containing $v^G_\emptyset(\sigma).$ By
definition $K$ is a compact open subgroup of $G$ and $\lambda$ is an
irreducible representation of $K,$ such that the functor 
$$V\mapsto Hom_G(c-i^G_K(\lambda),V)$$   
is an equivalence of categories from the block above to the category
of right $End_G(i^G_K(\lambda))$-modules.
There exists an unramified extension $F'/F$, such that the following holds (\cite{BK},\cite{V2}). Set
$G'=GL_k(F')$ and let
$I'\subset G'$ be the standard Iwahori subgroup. Then there is 
an algebra isomorphism \cite{BK} 7.6.19 
$$End_{G'}(i^G_{I'}({\bf 1}))
 \rightarrow End_G(i^G_K(\lambda)).$$
This isomorphism induces an equivalence between the block of unipotent
$G'$-represen\-tations
and the block of $G$-representations containing
$v^G_\emptyset(\sigma).$
Under this equivalence, the representations $v^G_I(\sigma)$ and
$v^G_{P_I}$ correspond to each other. This can be seen from the fact that the equivalence is
compatible with normalized induction \cite{BK} 7.6.21 and with twists
\cite{BK} 7.5.12. Thus, the statement follows from Corollary 2. \qed

\bigskip
\begin{flushleft}
Universit\"at Leipzig \\ 
Fakult\"at f\"ur Mathematik und Informatik \\ 
Mathematisches Institut \\ 
Augustusplatz 10/11 \\ 
D-04109 Leipzig \\ 
Germany 
\end{flushleft}

\end{document}